\documentclass[12pt,reqno]{amsart}
\usepackage[pagewise]{lineno}
\usepackage[most]{tcolorbox}
\tcbset{breakable,arc=0mm,colframe=white,colback=black!7.5!white}
\usepackage{graphicx}
\usepackage{amsmath,amsthm,amssymb,enumerate}
\usepackage{mathtools}
\usepackage{euscript,mathrsfs}
\usepackage{dsfont}
\usepackage{xcolor}
\usepackage{empheq}
\usepackage[left=2cm,right=2cm,top=3.5cm,bottom=3.5cm]{geometry}
\usepackage{color}
\catcode`\@=11 \@addtoreset{equation}{section}

\catcode`\@=12

\tcbset{arc=0mm,colframe=white,colback=black!10!white}

\newtheorem{Theorem}{Theorem}[section]
\newtheorem{Proposition}[Theorem]{Proposition}
\newtheorem{Assumptions}[Theorem]{Assumptions}
\newtheorem{Lemma}[Theorem]{Lemma}
\newtheorem{Corollary}[Theorem]{Corollary}

\theoremstyle{definition}
\newtheorem{Definition}[Theorem]{Definition}

\newtheorem{Remark}[Theorem]{Remark}

\newcommand{\bTheorem}[1]{
\begin{Theorem} \label{T#1} }
\newcommand{\eT}{\end{Theorem}}

\newcommand{\bProposition}[1]{
\begin{Proposition} \label{P#1}}
\newcommand{\eP}{\end{Proposition}}

\newcommand{\bLemma}[1]{
\begin{Lemma} \label{L#1} }
\newcommand{\eL}{\end{Lemma}}
\newcommand{\bCorollary}[1]{
\begin{Corollary} \label{C#1} }
\newcommand{\eC}{\end{Corollary}}

\renewcommand{\(}{\left(}
\renewcommand{\)}{\right)}

\allowdisplaybreaks

\renewcommand{\u}{\mathbf{u}}
\newcommand{\U}{\mathbf{U}}
\newcommand{\vv}{\mathbf{v}}
\newcommand{\ww}{\mathbf{w}}

\newcommand{\J}{\mathbf{J}}
\newcommand{\Jh}{\skew{2}\hat{\mathbf{J}}}
\newcommand{\bxi}{\boldsymbol{\xi}}
\newcommand{\dta}{\tau}
\newcommand{\ddta}{\,\mathrm{d}\tau}
\newcommand{\iQt}{\int_{\Omega_t}}

\newcommand{\Du}{\mathrm{D}\u}
\renewcommand{\div}[1]{\mathrm{div}\left({#1}\right)}
\newcommand{\di}[1]{\mathrm{div}\,{#1}}
\newcommand{\dib}[1]{\mathrm{div}\big({#1}\big)} 
\newcommand{\p}{\partial}
\newcommand{\pt}{\partial_t}
\newcommand{\na}{\nabla}
\newcommand{\R}{\mathbb{R}}

\newcommand{\dx}{\,\mathrm{d}x}

\newcommand{\dt}{\,\mathrm{d}t}

\DeclarePairedDelimiter{\norm}{\|}{\|}
\DeclarePairedDelimiter{\snorm}{|}{|}

\DeclareMathOperator{\spn}{span}

\newcommand\restr[2]{\ensuremath{\left.#1\right|_{#2}}}

\def\softd{{\leavevmode\setbox1=\hbox{d}%
          \hbox to 1.05\wd1{d\kern-0.4ex{\char039}\hss}}}

\date{}          

\title[Viscoelastic phase separation]{Existence and weak-strong uniqueness for global weak solutions for the viscoelastic phase separation model in three space dimensions}

\author{Aaron Brunk} 
\thanks{Numerical Mathematics, Department of Mathematics, Johannes Gutenberg University Mainz, Staudingerweg. 9, 55099 Mainz, Germany,
abrunk@uni-mainz.de}
\date{\today}

\begin{document}

\begingroup
\def\uppercasenonmath#1{} 
\let\MakeUppercase\relax 
\maketitle
\endgroup

\begin{abstract}
	\noindent  The aim of this work is to prove the global-in-time existence of weak solutions for a viscoelastic phase separation model in three space dimensions. To this end, we apply the relative energy concept provided by \cite{brunk2021relative}. We consider the case of regular polynomial-type potentials and positive mobilities, as well as the degenerate case with logarithmic potential and vanishing mobility.  
\end{abstract}


\section{Introduction}

The phase separation process of simple binary fluids is a central phenomenon in soft matter physics. A possibility to observe this is by quenching the mixture from a high-temperature regime to a low-temperature regime. For Newtonian
fluids, this phenomenon is well studied and the generic macroscopic model is the Cahn-Hilliard-Navier-Stokes system (model H) by Hohenberg and Halperin \cite{Hohenberg.1977}. \\

If considering a mixture where one constituent is much larger, for example, a polymer, the arising asymmetry plays an important role in the dynamics. In addition to the typical viscoelastic effects which are modelled by non-Newtonian rheology, the \emph{dynamic asymmetry} induces additional effects.
In this context, Tanaka \cite{Tanaka.200} introduced a mathematical model and coined the term viscoelastic phase separation exactly for mixtures that are governed by the effects of the above-mentioned dynamical asymmetric.  Unfortunately, this model was not consistent with the second
law of thermodynamics. In \cite{Zhou.2006} Zhou, Zhang, E re-derived a consistent model
using thermodynamic consistent methods, \cite{Grmela.1997,Grmela.1997b,Groot.2016}. A key ingredient of both models is an additional \emph{pressure} which is connected to the velocity difference of both multi-phase fluids, see \cite{Brunk2021} for more details on the closure relations. This pressure, later denoted $q$, was introduced by Tanaka \cite{Tanaka.200} accounting for inter-molecular interactions of the polymer-solvent mixture and is believed to be a viscoelastic phenomenon.  Note that both models are extensions of the model H.\\

In this work, we generalize our two-dimensional results for a viscoelastic phase separation model from \cite{Brunk.d,Brunk.c,brunk2021relative} to three-space dimensions. Furthermore, we will not consider the viscoelastic effects due to non-linear rheology but focus only on the impact of the additional \emph{pressure} on the system. We observe that the new pressure forms a so-called cross-diffusive subsystem with the Cahn-Hilliard-type equation. It is well-known that many cross-diffusion problems need a delicate treatment in the analysis, see e.g. \cite{Jungel.2016}. 
Note that the \emph{boundness-by-entropy} framework of Jüngel \cite{Jungel.2016} cannot be applied in our case, because it is so far only developed for second-order equations with a convex entropy. However, the Cahn-Hilliard equation is a fourth-order equation with the non-convex energy functional.\\

Similar models, typically without the cross-diffusion part, were already considered in the literature. For instance, the Cahn-Hilliard-Navier-Stokes model (model H) is discussed in \cite{Boyer.1999}. Here the author provides existence results in two and three space dimensions for the regular and the degenerate case and moreover higher regularity results for the regular case. We want to point out that so far there is no uniqueness result for the Cahn-Hilliard equation with a degenerate mobility (partially vanishing coefficient) and a logarithmic potential. For a broad review of the literature, we refer the reader to \cite{Cherfils.2011}.

We want to study a relative energy method that will finally result in a weak-strong uniqueness theorem. In the literature, one can find many results dealing with the compressible  Navier-Stokes equations \cite{Feireisl.2012,Feireisl.2016,Feireisl.2018}. The phase-field models with compressible Navier-Stokes equations were presented in \cite{Feireisl.2019,feireisl2021numerical}. A related model is the Navier-Stokes-Korteweg system, the corresponding results can be found in \cite{Giesselmann2017}. The relative energy and the corresponding weak-strong uniqueness results for the viscoelastic phase separation model are given in our recent work \cite{brunk2021relative}. There we gave a proof in two space dimensions and a conditional result in three space dimensions by requiring the existence of a suitable weak solution.\\

The structure of this work is as follows. In Section 2 we introduce the model in the regular case with polynomial potentials. Section 3 contains the theoretical framework and the notation which we employ. Section 4 will deal with the existence of weak solutions in the regular case. The focus in Section 5 is on the relative energy method and the weak-strong uniqueness result. In Section 6 we consider the degenerate case with logarithmic potentials, where we prove the existence of weak solutions and state a conditional weak-strong uniqueness result.

\section{Mathematical Model}
Let us start with introducing the viscoelastic phase separation model which will be concerned in this work. The model is based on a thermodynamic consistent framework where the total energy is given as
\begin{align}
E(\phi,q,\u)&=  E_{mix}(\phi) + E_{bulk}(q) + E_{kin}(\u) =\int_\Omega \(\frac{c_0}{2}\snorm*{\na\phi}^2 + F(\phi)\) + \int_\Omega \frac{1}{2}q^2  + \int_\Omega \frac{1}{2}\snorm*{\u}^2. \label{eq:free_energy}
\end{align}
We refer the reader to \cite{Brunk.d,Brunk2021} for a detailed derivation. Note that in contrast to \cite{Brunk.d} we do not consider the effects of the conformation tensor. Consequently, the corresponding system of partial differential equations reads

\begin{tcolorbox}
	\begin{align}
	\label{eq:full_model}
	\begin{split}
	\frac{\partial \phi}{\partial t} + \u\cdot\na\phi &= \dib{n^2(\phi)\nabla\mu} - \dib{n(\phi)\nabla\big(A(\phi)q\big)} \\
	\frac{\partial q   }{\partial t} + \u\cdot\na q &= -\frac{1}{\tau(\phi)}q + A(\phi)\Delta\big(A(\phi)q\big) - A(\phi)\dib{n(\phi)\nabla\mu} + \varepsilon_1\Delta q\\ 
	\frac{\partial\u   }{\partial t} + (\u\cdot\nabla)\u  &= \dib{\eta(\phi)\Du} -\nabla p  + \nabla\phi\mu \\
	\mu &= -c_0\Delta\phi + F^\prime(\phi),\quad \di\u = 0.
	\end{split}
	\end{align}
\end{tcolorbox}
System (\ref{eq:full_model}) is formulated on $(0,T)\times\Omega$, where, $\Omega\subset\mathbb{R}^3$ is sufficiently smooth. It is equipped with the following initial and boundary conditions
\begin{align}
\restr{(\phi,q,\u)}{t=0} = (\phi_0 , q_0, \u_0), \restr{\p_n\phi}{\p\Omega}=\restr{\p_n\mu}{\p\Omega}=\restr{\p_nq}{\p\Omega}=0,  \restr{\u}{\p\Omega} = \mathbf{0}. \label{eq:sysbc}
\end{align}

We proceed by imposing the following set of assumptions on the parameter functions.
\begin{Assumptions}
	\label{ass:paramreg}
	\hspace{1em}
	\begin{itemize}
		\item The functions $ n(s), \eta(s), \tau(s)$ are $C^0(\mathbb{R})$ and bounded by positive constants $e_1, n_2, \eta_1,\eta_2,\tau_1,\tau_2$ via $n_1\leq n(s) \leq n_2$, $\eta_1\leq \eta(s) \leq \eta_2$, $\tau_1\leq \tau(s) \leq \tau_2$. 
		\item We assume $A(s)\in C^1(\mathbb{R})$ and bounded by non-negative constants $A_1, A_2$ via $A_1\leq A(s) \leq A_2$ and $\norm{A^\prime}_{L^\infty(\mathbb{R})}\leq c$.
		\item We assume that $F\in C^2(\mathbb{R})$ with constants $c_{1,i}, c_{2,i} > 0,\; i=1,\ldots,2$ and $c_4\geq 0$ such that:
		\begin{align*}
			|F^{(i)}(x)| \leq c_{1,i}|x|^{p-i} + c_{2,i},\; i=0,1,2 \textnormal{ and } 4 \geq p\geq 2, \quad F(x) \geq -c_3, \quad F^{\prime\prime}(x) \geq - c_4. 
		\end{align*}
        \item The constants $\varepsilon_1, c_0$ are positive.
	\end{itemize}
\end{Assumptions}

For $\varepsilon_1=0$ the above model reduces to the simplified model of Zhou et al. \cite{Zhou.2006}. In this case, the necessary a priori estimates for the chemical potential $\mu$ cannot be obtained and the existence of weak solutions remains unclear. The additional term $\varepsilon_1\Delta q$ has also a physical motivation. The variable $q$ models an additional pressure and the governing equation is the scalar version for matrix-valued viscoelastic stresses, here Oldroyd-B type, see \cite{Brunk2021,Zhou.2006}. In the context of viscoelastic stresses, especially dilute polymer solutions, the Laplacian models center-of-mass diffusion of the polymer chains, see \cite{Barrett,MichalBathory.2020}

\section{Theoretical framework}
In this section, we introduce the notation and theoretical framework for our upcoming analysis of (\ref{eq:full_model}).  Further, we denote by $\Omega_T:=\Omega\times(0,T)$ and $\Omega_t:=\Omega\times(0,t)$ the full space-time cylinder and the intermediate space-time cylinder, respectively. The standard Lebesgue spaces are denoted by $L^p(\Omega)$ and their norm is denoted by $\norm*{\cdot}_p$. We use the standard notation for the Sobolev spaces and introduce the notation $
V:=H^1_{0,\text{div}}(\Omega)^3,  H:=L^2_{\text{div}}(\Omega)^3$. As usual, these spaces are obtained as closure with respect to the $L^2(\Omega)^3$, $H^1(\Omega)^3$ norm of the space of $C_{0,\infty}(\bar\Omega)^3$ function with zero divergence, respectively.

The space $V$ is equipped with the norm $\norm*{\cdot}_V:=\norm*{\na\cdot}_2$. We denote the dual space of $H^1(\Omega), W^{1,p}(\Omega), V$ by $(H^{1}(\Omega))^*, (W^{1,p}(\Omega))^*, V^*$ respectively and we denote the norms of the corresponding Bochner space $L^p(0,T;L^q(\Omega))$ by $\norm*{\cdot}_{L^p(L^q)}.$ 

Next, we state an interpolation lemma which we use in the analysis.
\begin{Lemma}[\cite{Grun.1995}]
	\label{lema:interpol}
Let $d$ be the space dimension and $v\in L^\infty(0,T;L^p(\Omega))\cap L^p(0,T;W^{1,p}(\Omega))$. Then, $v\in L^r(0,T;L^q(\Omega))$ and the estimate
\begin{equation*}
\norm*{v}_{L^r(L^q)}\leq c\(\norm*{v}_{L^\infty(L^p)} + \norm*{v}_{L^p(W^{1,p})}\)
\end{equation*}
holds. The constant c depends on the $d,p,|\Omega|,T$ and $\partial\Omega$. Moreover, $r,p,q$ and d have to satisfy the relations
\begin{equation*}
\frac{1}{r} = \frac{d}{p^2}-\frac{d}{pq}
\end{equation*} 
\begin{align*}
q\in\left[p,\frac{dp}{d-p}\right] \text{ and } r\in[p,\infty) \text{ if } 1 < p < d, \qquad q\in[p,\infty) \text{ and } r\in\left(\frac{p^2}{d},\infty\right] \text{ if } 1 < d \leq p.
\end{align*}
\end{Lemma}

The above lemma yields in three space dimensions for $p=2$ bounds in the space $L^4(0,T;L^3(\Omega))$ and $L^{10/3}(\Omega_T)$.

 \section{Weak solution and Existence}
 In this section, we introduce the notion of a weak solution for (\ref{eq:full_model}). Further, we will state the corresponding existence result.
\begin{Definition}
	\label{defn:weak_sol_full}
	Given the initial data $\left(\phi_0, q_0, \u_0\right) \in \big[H^1(\Omega)\times L^2(\Omega)\times H\big]$.
	The quadruple $(\phi,q,\mu,\u)$ is called a weak solution of \eqref{eq:full_model} if
	\begin{align*}
	&\phi \in L^{\infty}(0,T;H^1(\Omega))\cap L^2(0,T;H^2(\Omega))\cap H^1(0,T;(H^{1}(\Omega))^*),\\
	&q\in L^\infty(0,T;L^2(\Omega))\cap L^2(0,T;H^1(\Omega))\cap W^{1,4/3}(0,T;(H^{1}(\Omega))^*), \\
	& \u\in L^{\infty}(0,T;H)\cap L^2(0,T;V)\cap W^{1,4/3}(0,T;V^*),\\
	&A(\phi)q,\mu \in L^{5/3}(0,T,W^{1,5/3}(\Omega)),\\
	&n(\phi)\na\mu - \na\big(A(\phi)q\big) \in L^2(0,T;L^2(\Omega)),
	\end{align*}
	and  
	\begin{align}
	&\int_\Omega\pt\phi\psi \dx + \int_\Omega (\u\cdot\na\phi)\psi\dx + \int_\Omega n(\phi)\Big[n(\phi)\nabla\mu-\nabla(A(\phi)q)\Big]\cdot\na\psi \dx = 0 \label{eq:weak_sol_full}\\
	&\int_\Omega\pt q\zeta \dx  + \int_\Omega (\u\cdot\na q)\zeta\dx+ \int_\Omega\frac{q\zeta}{\tau(\phi)}\dx + \varepsilon_1\int_\Omega \na q \cdot\na\zeta \dx\nonumber\\
	&\hspace{12em}- \int_\Omega \Big[n(\phi)\nabla\mu - \nabla(A(\phi)q)\Big]\cdot\na(A(\phi)\zeta) \dx= 0 \nonumber\\ 
	&\int_\Omega\mu\xi \dx - c_0\int_\Omega\nabla\phi\cdot\nabla\xi \dx - \int_\Omega F^\prime(\phi)\xi \dx = 0 \nonumber\\
	&\int_\Omega\pt\u\cdot\vv \dx + \int_\Omega(\u\cdot\nabla)\u\cdot\vv \dx - \int_\Omega\eta(\phi)\Du:\mathrm{D}\vv \dx  + \int_\Omega\phi\nabla\mu\cdot\vv \dx = 0. \nonumber
	\end{align}
	The above system holds for any test function $(\psi,\zeta,\xi,\vv)\in [H^1(\Omega) \times H^1(\Omega)  \times H^1(\Omega) \times V]$ and almost every $t\in(0,T)$. Furthermore, it satisfies the initial conditions given above, i.e. $\left(\phi(0), q(0), \u(0) \right) = \left(\phi_0, q_0, \u_0 \right)$.
\end{Definition}

\begin{tcolorbox}
	\begin{Theorem}
		\label{theo:ex_full}
		Let Assumptions \ref{ass:paramreg} hold. For given initial data $(\phi_0, q_0, \u_0)\in [H^1(\Omega)\times L^2(\Omega)\times H]$ and time $T>0$ there exists a global-in-time weak solution of (\ref{eq:full_model}) in the sense of Definition \ref{defn:weak_sol_full} which satisfies for a.a. $t\in(0,T)$ the energy inequality 
		\begin{align}
		& \(\int_\Omega \frac{c_0}{2}\snorm*{\na\phi(t)}^2 + F(\phi(t)) + \frac{1}{2}\snorm*{q(t)}^2 + \frac{1}{2}\snorm*{\u(t)}^2  \dx \) \label{eq:fullie1_def} \\
		&+ \iQt \snorm*{n(\phi)\na\mu - \na\big(A(\phi)q\big) }^2 + \frac{1}{\tau(\phi)}|q|^2 + \varepsilon_1\snorm*{\na q}^2 +  \eta(\phi)\snorm*{\Du}^2\dx\ddta\nonumber \\
		&\leq  \(\int_\Omega \frac{c_0}{2}\snorm*{\na\phi(0)}^2 + F(\phi(0))+ \frac{1}{2}\snorm*{q(0)}^2 + \frac{1}{2}\snorm*{\u(0)}^2 \dx \). \nonumber
		\end{align}	
		
	\end{Theorem}	
\end{tcolorbox}

\begin{Remark}
	The proof is similar to \cite[Proof of Theorem 4.1]{Brunk.d}. First a suitable Galerkin approximation of (\ref{eq:weak_sol_full}) is introduced. By means of energy methods, we derive the necessary a priori estimates. Finally, we pass to the limit in the Galerkin approximation and obtain the weak formulation (\ref{eq:weak_sol_full}).
\end{Remark}

 \subsection{Approximation and Estimates}

In order to prove existence of a weak solution (\ref{eq:weak_sol_full}) of model (\ref{eq:full_model}) we consider a Galerkin approximation of problem (\ref{eq:weak_sol_full}).
Let $\psi_j,  j=1,\ldots,\infty,$ be smooth basis functions such that
\begin{align*}
H^1(\Omega) = \overline{\spn\{\psi_j\}_{j=1}^{\infty}},\quad  V=\overline{\spn\{\vv_j\}_{j=1}^{\infty}}.
\end{align*}
 Furthermore, $\psi_j$ are eigenfunctions of the negative Laplacian $-\Delta$ subjected to Neumann boundary conditions and $\vv_j$ are divergence-free functions subjected to Dirichlet boundary conditions. For simplicity, they are chosen to be orthonormal in $L^2(\Omega), H$ and orthogonal in $H^1(\Omega), V$. We introduce the finite-dimensional subspaces $W_m:=\spn\{\psi_1,\ldots,\psi_m\}, V_m:=\spn\{\vv_1,\ldots,\vv_m\}$ and the corresponding orthogonal projection by $\mathbf{P}_{W_m}, \mathbf{P}_{V_m}$, respectively. We define the $m$-th Galerkin approximation via
\begin{align*}
\phi_m(x,t) &= \sum_{j=1}^m \lambda_{jm}(t)\psi_j(x), \hspace{1em} \mu_m(x,t) = \sum_{j=1}^m \theta_{jm}(t)\psi_j(x),\\
 q_m(x,t) &= \sum_{j=1}^m \zeta_{jm}(t)\psi_j(x),\hspace{1em} \u_m(x,t) = \sum_{j=1}^m g_{jm}(t)\vv_j(x), \\
\phi_m(0) &= \mathbf{P}_{W_m}(\phi_{0}),\hspace{2em} q_m(0) = \mathbf{P}_{W_m}(q_{0}),\hspace{2em} \u_m(0) = \mathbf{P}_{V_m}(\u_{0}).
\end{align*}

The Galerkin approximations satisfies $(\ref{eq:weak_sol_full})_m$, i.e.,
\begin{align}
	&\int_\Omega\pt\phi_m\psi_m \dx + \int_\Omega (\u_m\cdot\na\phi_m)\psi_m\dx + \int_\Omega n(\phi_m)\Big[n(\phi_m)\nabla\mu_m-\nabla(A(\phi_m)q_m)\Big]\cdot\na\psi_m \dx = 0 \nonumber\\
	&\int_\Omega\pt q_m\zeta_m \dx  + \int_\Omega (\u_m\cdot\na q_m)\zeta_m\dx+ \int_\Omega\frac{q_m\zeta_m}{\tau(\phi_m)}\dx + \varepsilon_1\int_\Omega \na q_m \cdot\na\zeta_m \dx\nonumber\\
	&\hspace{12em}- \int_\Omega \Big[n(\phi_m)\nabla\mu_m - \nabla(A(\phi_m)q_m)\Big]\cdot\na(A(\phi_m)\zeta_m) \dx= 0 \nonumber \\ 
&\int_\Omega\mu_m\xi_m \dx - c_0\int_\Omega\nabla\phi_m\cdot\nabla\xi_m \dx - \int_\Omega F^\prime(\phi_m)\xi_m \dx = 0 \label{eq:weak_sol_fullm}\\
	&\int_\Omega\pt\u_m\cdot\vv_m \dx + \int_\Omega(\u_m\cdot\nabla)\u_m\cdot\vv_m \dx - \int_\Omega\eta(\phi_m)\Du_m:\mathrm{D}\vv_m \dx  + \int_\Omega\phi_m\nabla\mu_m\cdot\vv_m \dx = 0. \nonumber
	\end{align}
for all $(\psi_m, \zeta_m,\xi_m,\vv_m) \in W_m\times W_m\times W_m\times V_m$.

By standard techniques from ordinary differential equations, the solutions exist up to time $T_m$. 

\subsection{A priori estimates}
We reproducing the discrete version of energy inequality (\ref{eq:fullie1_def}) by inserting $\psi_m=\mu_m, \zeta_m=q_m, \xi_m=\partial_t \phi_m, \vv_m=\u_m$ into \eqref{eq:weak_sol_fullm}. The computations can be found in  \cite[Theorem 3.2] {Brunk.c} and yields
\begin{align}
E_m(t):=&\left(\int_\Omega \frac{c_0}{2}\snorm{\nabla\phi_m(t)}^2 + F(\phi_m(t)) + \frac{1}{2}|q_m(t)|^2 + \frac{1}{2}\snorm*{\u_m(t)}^2 \dx\right)  \label{eq:energy_pqm_disc}\\
& +\iQt \snorm*{n(\phi_m)\na\mu_m - \na\big(A(\phi_m)q_m\big)}^{2} +\varepsilon_1\snorm*{\nabla q_m}^2  + \frac{1}{\tau(\phi_m)}q_m^2 + \eta(\phi_m)\snorm*{\Du_m}^2\dx\ddta  \nonumber  \\
\leq&  \left(\int_\Omega \frac{c_0}{2}\snorm{\nabla\phi_m(0)}^2 + F(\phi_m(0)) + \frac{1}{2}|q_m(0)|^2 + \frac{1}{2}\snorm*{\u_m(0)}^2\dx\right). \nonumber
\end{align}
An application of the Gronwall Lemma to \eqref{eq:energy_pqm_disc} implies that $T_m=T$ for every $m$ and we obtain  a priori estimates independent of $m$ in the following spaces
\begin{align}
&\phi_m\in L^\infty(0,T;H^1(\Omega)), \quad q_m\in L^\infty(0,T;L^2(\Omega))\cap L^2(0,T;H^1(\Omega)),\label{eq:pqmapriori}\\ 
&\u_m \in L^\infty(0,T;L^2(\Omega))\cap L^2(0,T;V), \quad  n(\phi_m)\na\mu_m-\na\big(A(\phi_m)q_m\big) \in L^2(0,T;L^2(\Omega)). \nonumber
\end{align}

For the above result, control of the mean value of $\phi_m$ follows from mass-conservation, i.e. $\psi_m=1$. One can observe that $\int_\Omega \phi_m(t) = \int_\Omega \phi_m(0)$ for all $t\in(0,T)$. The complete estimate in $H^1$ follows from Poincare inequality.\\

\noindent\textbf{Refined estimates:} \\
For the limiting process, we need a priori estimates for $\mu_m$. Using \eqref{eq:pqmapriori} we can deduce that $\nabla(A(\phi_m)q_m)$ is bounded in $L^2(0,T;L^{3/2}(\Omega))$. This implies that $n(\phi_m)\nabla\mu_m$ and hence $\nabla\mu_m$ is also bounded in $L^2(0,T;L^{3/2}(\Omega))$.

In order to obtain a bound for $\mu_m$ we consider its mean value, i.e., by testing \eqref{eq:weak_sol_fullm} with $\xi=1$ to obtain
\begin{align}
(\mu_m)_\Omega = (F^{\prime}(\phi_m))_\Omega \nonumber
\Longrightarrow \int_0^T \snorm*{(\mu_m)_\Omega}^q \leq 
\int_0^T \norm*{\phi_m}_{p-1}^{q(p-1)}  \leq C\norm*{\phi_m}^{(p-1)q}_{L^\infty(L^{p-1})}\leq C \label{eq:pqmmumean} 
\end{align}
for some $q\in[1,\infty)$. The above bound follows, since $p\leq 4$, cf. Assumption \ref{ass:paramreg}, and $\phi_m$ is bounded in $L^\infty(0,T;H^1(\Omega))$.
This implies that $(\mu_m)_\Omega$ is bounded in $L^q(0,T)$ for all $q\in [1,\infty)$. By applying the $L^p$-version of the Poincaré inequality we get
\begin{equation}
\norm*{\mu_m-(\mu_m)_\Omega}_{L^2(L^{3/2})} \leq c\norm*{\na\mu_m}_{L^2(L^{3/2})} \leq C. \label{eq:pqmmuconv}
\end{equation}
 Hence we obtain $\mu_m\in L^{2}(0,T;W^{1,3/2}(\Omega))$ and by the Sobolev embedding $\mu_m\in L^{2}(0,T;L^{3}(\Omega))$. 
 
With these estimates we choose $\xi_m=-\Delta\phi_m$, integrate over time $(0,T)$ to derive
\begin{align*}
   \int_0^T\norm{\Delta\phi_m}_{2}^2 &\leq C\int_0^T\norm{\mu_m}_2^2 + \norm{F'(\phi_m)}_2^2 \dt. \\
   & \leq C + C\int_0^T\norm{\phi_m}_{6}^6 \dt \leq C.
\end{align*}
Thus, we deduce by norm equivalences that $\phi_m\in L^2(0,T;H^2(\Omega))$.
To summarize we obtain refined a priori bounds independent of $m$ in the following spaces
\begin{align}
  A(\phi_m)q_m , \mu_m  \in L^2(0,T;W^{1,3/2}(\Omega)), \quad \phi_m\in L^2(0,T;H^2(\Omega)). \label{eq:apriorinew2}
\end{align}

\noindent\textbf{Bootstrapping:} \\
We will now exploit the refined a priori bounds to increase the regularity. First, we reconstruct some crucial estimates from the cross-diffusive part of the energy dissipation \eqref{eq:energy_pqm_disc}, i.e. from $\snorm*{n(\phi_m)\na\mu_m-\na\big(A(\phi_m)q_m \big)}^2$. By Lemma \ref{lema:interpol} we find $q_m,\na\phi_m\in L^{10/3}(\Omega_T)\cap L^4(0,T;L^3(\Omega))\cap L^{8/3}(0,T;L^4(\Omega))$. We calculate for $p\in [1,\infty)$
\begin{align*}
\int_{\Omega_T} |\na\big(A(\phi_m)q_m\big)|^p &\leq c\int_{\Omega_T}\norm*{\na\phi_m q_m}^p + \norm*{\na q_m}^p \leq c\norm*{\na\phi_mq_m}^p_{L^p(\Omega_T)} + c\norm*{\na q_m}^p_{L^p(\Omega_T)}.
\end{align*}
Using the obtained regularity (\ref{eq:pqmapriori}) we observe that the second term is bounded for $p\leq 2$. Lemma \ref{lema:interpol} and the generalized Hölder inequality imply that the first term is bounded for $p\leq 5/3$. Hence, we find
\begin{align}
\na\big(A(\phi_m)q_m \big)\in L^{5/3}(\Omega_T). \label{eq:pqmAqreg}
\end{align}
By virtue of $\eqref{eq:pqmapriori}_2$ together with \eqref{eq:pqmAqreg} we obtain $\na\mu_m \in L^{5/3}(\Omega_T)$. As $(\mu_m)_\Omega$ is also bounded in $L^{5/3}(0,T)$ we derive the new a priori bounds in the space
\begin{equation}
A(\phi_m)q_m, \mu_m \in L^{5/3}(0,T;W^{1,5/3}(\Omega)).    \label{eq:add_aprior}
\end{equation}

\subsection{Time derivative \& Compact embedding}

Using orthogonality of the eigenfunctions one can deduce that
\begin{align}
   \norm{\pt\phi_m}_{(H^1)^*} = \sup_{\psi\in H^1} \frac{\int_\Omega \pt\phi_m \psi}{\norm{\psi}_{H^1}} = \sup_{\psi_m\in W_m} \frac{\int_\Omega \pt\phi_m \psi_m}{\norm{\psi_m}_{H^1}}, \label{eq:tdre1}\\
   \norm{\pt q_m}_{(H^1)^*} = \sup_{\psi\in H^1} \frac{\int_\Omega \pt q_m \psi}{\norm{\psi}_{H^1}} = \sup_{\psi_m\in W_m} \frac{\int_\Omega \pt q_m \psi_m}{\norm{\psi_m}_{H^1}}, \label{eq:tdre2}\\
   \norm{\pt\u_m}_{V^*} = \sup_{\vv\in V} \frac{\int_\Omega \pt\u_m \vv}{\norm{\vv}_{H^1}} = \sup_{\vv_m\in V_m} \frac{\int_\Omega \pt\u_m \vv_m}{\norm{\vv_m}_{V}}.\label{eq:tdre3}
\end{align}

To obtain a priori bounds on the time derivative, we will further estimate the above norms using the weak formulation. Considering \eqref{eq:tdre1} we find
\begin{align}
\int_0^T \norm*{\pt \phi_m}^2_{(H^{1})^*}\dt &\leq c\int_0^T \norm*{n(\phi_m)\na\mu_m-\na\big(A(\phi_m)q_m\big)}_{2}^2 + \norm*{\u_m}_3^2\norm*{\na\phi_m}_2^2 \dt \leq C_0.
\label{eq:pqm_TDphi}
\end{align}
We observe that our a priori bounds \eqref{eq:pqmapriori} together with the interpolation Lemma \ref{lema:interpol} imply  that the above integral is bounded,  which yields $\partial_t\phi_m\in L^2(0,T;(H^{1}(\Omega))^*)$.

In view of \eqref{eq:tdre2} we estimate
\begin{align*}
\int_0^T \norm*{\pt q_m}^{4/3}_{(H^{1})^*} \dt\leq&\; c\int_0^T \norm*{q_m}_2^{4/3} + \varepsilon_1\norm*{\na q_m}_2^{4/3} + \norm*{n(\phi_m)\na\mu_m-\na\big(A(\phi_m)q_m\big)}_{2}^{4/3} \dt \label{eq:pqm_TDq}\\
&+ c\int_0^T \norm*{n(\phi_m)\na\mu_m-\na\big(A(\phi_m)q_m\big)}_{2}^{4/3}\norm*{\na\phi_m}_{3}^{4/3} + \norm*{\u_m}_{3}^{4/3}\norm*{\na q_m}_2^{4/3} \dt \nonumber \\ 
\leq& \;\norm*{q_m}_{L^{4/3}(L^2)}^{4/3} + \varepsilon_1\norm*{\na q_m}_{L^{4/3}(L^2)}^{4/3} + \norm*{n(\phi_m)\na\mu_m-\na\big(A(\phi_m)q_m\big)}_{L^{2}(L^{2})}^2 \nonumber\\
&+ \norm*{n(\phi_m)\na\mu_m-\na\big(A(\phi_m)q_m\big)}_{L^{2}(L^{2})}^2 + \norm*{\na\phi_m}_{L^{4}(L^{3})}^4 + \norm*{\u_m}_{L^4(L^3)}^{4} + \norm*{\na q_m}_{L^{2}(L^{2})}^2\nonumber \\
 \leq &\;C_0. \nonumber
\end{align*}
The above calculations together with the regularities in \eqref{eq:pqmapriori} show that the time derivative $\pt q_m$ remains bounded in $L^{4/3}(0,T;(H^{1}(\Omega))^*)$.

Before estimating the last time derivative we estimate the following term
\begin{align}
\int_0^T \int_\Omega \phi_m\nabla\mu_m\cdot\vv_m  \dx\dt &\leq \int_0^T \norm*{\phi_m}_{6}\norm*{\na\mu_m}_{3/2}\norm*{\vv_m}_6 \dt, \\
\int_0^T \norm*{\phi_m}_{4/3}^{6}\norm*{\na\mu_m}_{3/2}^{4/3}\dt &\leq c\norm*{\phi_m}^{4/3}_{L^\infty(H^1)}\norm*{\na\mu_m}_{L^{4/3}(L^{3/2})}^{4/3} \leq C_0.
\end{align}
Finally, we consider \eqref{eq:tdre3} and find
\begin{align*}
\int_0^T \norm{\pt\u_m}_{V^*}^{4/3} \dt \leq C\int_0^T\norm{\mathrm{D}\u_m}_2^{4/3} + \norm{\u_m}_4^{8/3} \dt + \norm*{\phi_m}^{4/3}_{L^\infty(H^1)}\norm*{\na\mu_m}_{L^{4/3}(L^{3/2})}^{4/3}.  
\end{align*}
Again using a priori estimates \eqref{eq:pqmapriori}, \eqref{eq:add_aprior} together with the interpolation Lemma \ref{lema:interpol} we can see that the time derivative $\pt \u_m$ remains bounded in $L^{4/3}(0,T;V^*)$.\\

\noindent\textbf{Compact embeddings:}\\
Due to the Lemma of Banach-Alaoglu and the Lemma of Aubin-Lions, we obtain the following convergences after suitable extraction of subsequences if necessary.
\begin{align}
\phi_m &\rightharpoonup^{\star} \phi \in L^\infty(0,T;H^1(\Omega)),& q_m &\rightharpoonup^{\star} q \in L^\infty(0,T;L^2(\Omega)), \label{eq:pqmconvtable}\\
\phi_m &\rightharpoonup \phi \in L^2(0,T;H^2(\Omega))\cap L^{10/3}(0,T;W^{1,10/3}(\Omega)),& q_m &\rightharpoonup q \in L^2(0,T;H^1(\Omega))\cap L^{10/3}(\Omega_T),\nonumber \\
\phi_m &\rightarrow \phi \in L^{p}(0,T;W^{1,p}(\Omega)) \text{ for } p<\frac{10}{3},& q_m &\rightarrow q \in L^{p}(0,T;L^p(\Omega)) \text{ for } p<\frac{10}{3} ,\nonumber\\
\phi_m &\rightarrow \phi \in  L^2(0,T;W^{1,p}(\Omega)) \text{ for } p<6,& q_m &\rightarrow q \in  L^2(0,T;L^p(\Omega)) \text{ for } p<6,\nonumber \\
\pt\phi_m &\rightharpoonup \pt\phi \in L^{2}(0,T;(H^{1}(\Omega))^*),& \pt q_m &\rightharpoonup \pt q \in L^{4/3}(0,T;(H^{1}(\Omega))^*),\nonumber
\end{align}
\begin{align}
\u_m &\rightharpoonup^{\star} \u \in L^\infty(0,T;L^2(\Omega)), & \label{eq:uconvtable}\\
\u_m &\rightharpoonup \u\in L^2(0,T;V)\cap L^{10/3}(\Omega_T),&\nonumber \\
\u_m &\rightarrow \u \in  L^2(0,T;L^{p}(\Omega))\cap L^{q}(\Omega_T) \text{ for } p<6 \text{ and } q <\frac{10}{3},&\nonumber \\
\pt\u_m &\rightharpoonup \pt\u \in L^{4/3}(0,T;V^*).\nonumber&
\end{align}
Furthermore, due to the strong convergence $\phi_m, q_m$ and $\u_m$ converge almost everywhere in $\Omega_T$.\\

\noindent\textbf{Nonlinear limits:}\\
Before passing to the limit, we consider the weak limit of $\na\big(A(\phi_m)q_m)\big)$. To this end we calculate for a $\zeta\in L^q(\Omega_T)$ with $q\in[1,6]$ the integral
\begin{align*}
\int_0^T\int_\Omega A(\phi_m)\na q_m\zeta + A^\prime(\phi_m)\na\phi_mq_m\zeta.
\end{align*}
Since $A(\phi_m)\zeta$ converges strongly to $A(\phi)\zeta$ in $L^q(\Omega_T)$, cf. (\ref{eq:pqmconvtable}), we can pass to the limit for the first term. For the second term we see by (\ref{eq:pqmconvtable}) that the sequences $\{\na\phi_m\}$ and $\{q_m\}$ are bounded in $L^{10/3}(\Omega_T)$. Consequently, their product is weakly convergent in $L^{5/3}(\Omega_T)$. Further, the product $\na\phi_mq_m$ converges strongly in at least $L^1(\Omega_T)$ to $\na\phi q$, cf. \eqref{eq:pqmconvtable}. From the uniqueness of the weak limit, we conclude that $\na\phi_mq_m$ converges weakly to $\na\phi q$ in $L^{5/3}(\Omega_T)$. This yields the convergence, if $q\geq  5/2$. Finally, this implies
\begin{equation}
\na\big(A(\phi_m)q_m\big) \rightharpoonup \na\big(A(\phi)q\big) \in L^{5/3}(\Omega_T). \label{eq:pqmaphiconv53}
\end{equation}

In order to pass to the limit in the Cahn-Hilliard equation $(\ref{eq:weak_sol_full})_1$ and the energy inequality (\ref{eq:energy_pqm_disc}) we need the weak convergence of $n(\phi_m)\na\mu_m-\na\big(A(\phi_m)q_m \big)$ in $L^2(\Omega_T)$. In contrast to the two-dimensional case in \cite{Brunk.d} this does not follow directly from (\ref{eq:pqmconvtable}) and (\ref{eq:pqmaphiconv53}).  A priori we only know that the desired limit is bounded in $L^{5/3}(\Omega_T)$, cf. (\ref{eq:pqmaphiconv53}).
This implies the following weak convergences
\begin{align}
n(\phi_m)\na\mu_m - \na\big(A(\phi_m)q_m \big) &\rightharpoonup \ww \in L^2(\Omega_T), \nonumber\\
n(\phi_m)\na\mu_m - \na\big(A(\phi_m)q_m \big) &\rightharpoonup n(\phi)\na\mu - \na\big(A(\phi)q \big) \in L^{5/3}(\Omega_T).\label{eq:pqmcrossdifflimit}
\end{align}
By embedding $\ww$ into $L^{5/3}(\Omega_T)^d$ we obtain $\ww=n(\phi)\na\mu - \na\big(A(\phi)q \big) + \boldsymbol{\theta} \in L^2(\Omega_T)$ with 
\begin{equation}
\int_{\Omega_T} \boldsymbol{\theta}\cdot\boldsymbol{\zeta} = 0 \text{ for all } \boldsymbol{\zeta}\in L^{5/2}(\Omega_T)^d. \label{eq:pqmthetadefect}
\end{equation}
As $\snorm*{\Omega_T}<\infty$ the characteristic functions $\chi_E(x,t)$ are dense all $L^{p}(\Omega_T)$, $p< \infty$, for every measurable set $E\subset \Omega_T$. Inserting $\boldsymbol{\zeta}_j=\chi_E(x,t)$  in \eqref{eq:pqmthetadefect} yields
\begin{equation*}
\int_{E} \boldsymbol{\theta}_j= 0 \text{ for all measurable } E\subset\Omega_T \text{ and every } j=1,\ldots,3.
\end{equation*}
This implies $\boldsymbol{\theta}=\mathbf{0}$ a.e. in $\Omega_T$ and we find
\begin{equation}
n(\phi_m)\na\mu_m - \na\big(A(\phi_m)q_m \big) \rightharpoonup n(\phi)\na\mu - \na\big(A(\phi)q \big) \in L^{2}(\Omega_T).\label{eq:pqmcrossreallim}
\end{equation}

\subsection{Limit passing}
In the following we pass to the limit $m\to \infty$ in the main nonlinearities of $\eqref{eq:weak_sol_full}_m$. This is done by multiplying $\eqref{eq:weak_sol_full}_m$ by a time-dependent test function $\varphi \in L^\infty(0,T)$. Note that technically speaking, the test function for the discrete solution is taken from $W_m, V_m$ while the desired spaces for the weak solution are $H^1(\Omega), V$. By suitable insertion of a zero, this can be reduced to the limit in the test functions and the limit in the approximate solutions. The limit in the test functions is typically omitted since, they converge strongly in $H^1(\Omega), V$ by definition. Hence, we focus on the limit in the approximations.
In the Cahn-Hilliard equation $(\ref{eq:weak_sol_full})_{1,m}$ we only have to consider
\begin{align}
P_{1,m}:=\int_0^T\int_\Omega \Big[n(\phi_m)\big(n(\phi_m)\na\mu_m - \na\big(A(\phi_m)q_m \big)\big) - n(\phi)\big(n(\phi)\na\mu + \na\big(A(\phi)q \big) \big)\Big]\cdot\na\psi\varphi(t) \dx\dt \label{eq:pqmchlimit}
\end{align}
Since $n(\phi_m)$ to $n(\phi)$ converges a.e. in $\Omega_T$ we can apply the weak convergence of $n(\phi_m)\na\mu_m - \na\big(A(\phi_m)q_m \big)$ in $L^2(\Omega_T)$, cf. (\ref{eq:pqmcrossreallim}) to obtain $P_{1,m}\to 0$ as $m\to\infty$.

The main nonlinearity of $(\ref{eq:weak_sol_full})_{2,m}$ is given by
\begin{align*}
&-\int_0^T\int_\Omega \(n(\phi_m)\na\mu_m - \na\big(A(\phi_m)q_m\big) \)\na\big(A(\phi_m)\zeta\big)\varphi \\
= &\int_0^T\int_\Omega \(n(\phi_m)\na\mu_m - \na\big(A(\phi_m)q_m\big) \) \(A(\phi_m)\na\zeta + A^\prime(\phi_m)\na\phi_m\zeta \)\varphi := P_{2,m} + P_{3,m}. 
\end{align*}

Using the weak convergence of $n(\phi_m)\na\mu_m - \na\big(A(\phi_m)q_m\big)$ in $L^{2}(\Omega_T)$ we only have to show that $A(\phi_m)\na\zeta\varphi$ and $A^\prime(\phi_m)\na\phi_m\zeta\varphi$ converge strongly in $\in L^{2}(\Omega_T)$. As $A(\phi_m),A^\prime(\phi_m)$ converges to $A(\phi),A^\prime(\phi)$ a.e. in $\Omega_T$ we can conclude that $A(\phi_m)\na\zeta\varphi$ converges strongly to $A(\phi)\na\zeta\varphi$ in $L^2(\Omega_T)$. This implies $P_{2,m}\to P_2$ as $m\to\infty$. Similarly, by virtue of \eqref{eq:pqmconvtable} we observe that $A^\prime(\phi_m)\na\phi_m\zeta\varphi$ converges strongly to $A^\prime(\phi)\na\phi\zeta\varphi$ in $L^{2}(\Omega_T)$. Hence, we found $P_{3,m}\to P_3$ as $m\to\infty$. 
All other terms of $\eqref{eq:weak_sol_full}_{2,m}$ as well as the terms in $\eqref{eq:weak_sol_full}_{3,m}$ are treated similarly as in \cite{Brunk.d} and are omitted here.

In the last term, we consider the capillary stress tensor in the Navier-Stokes equation $\eqref{eq:weak_sol_full}_4$ given by
\begin{align}
&\int_0^T\int_\Omega\(\phi_m\na\mu_m - \phi\na\mu\)\vv\varphi \dx\dt = \int_0^T\int_\Omega\( (\phi_m-\phi)\na\mu + \na(\mu_m-\mu)\phi_m \)\vv\varphi \dx\dt \nonumber\\
&\leq \int_0^T \norm*{\phi_m- \phi}_{30/7}\norm*{\na\mu}_{5/3}\norm*{\vv}_6\norm*{\varphi}_\infty \dt + \int_0^T\int_\Omega \na(\mu_m-\mu)\phi_m\vv\varphi \dx\dt \nonumber\\
&\leq \norm*{\phi_m-\phi}_{L^{5/2}(H^1)}\norm*{\na\mu}_{L^{5/3}(L^{5/3})}\norm*{\vv}_6\norm*{\varphi}_\infty + \int_0^T\int_\Omega \na(\mu_m-\mu)\phi_m\vv\varphi \dx\dt. \label{eq:full_phistress}
\end{align}
The first integral of (\ref{eq:full_phistress}) tends to zero due to the strong convergence of $\phi_m$ in $L^3(0,T;H^1(\Omega))$, cf. \eqref{eq:pqmconvtable}. The second integral of \eqref{eq:full_phistress} goes to zero by weak convergence, since $\phi_m\vv\varphi$ converge strongly to $\phi\vv\varphi$ in $L^{5/2}(\Omega_T)$, cf. \eqref{eq:uconvtable}. Note that all other terms can be treated by similar arguments as in \cite[Section 8]{Brunk.d} and hence we omit them here.\\[0.5em]

\subsection{Energy limit}
In order to pass to the limit in the energy equality \eqref{eq:energy_pqm_disc} we recall that for a suitable weakly/weakly-$^\star$ convergent sequence $\{g_m\}$ we have
\begin{align}
\norm*{g(t)}_2 \leq \norm*{g}_{L^\infty(0,t;L^2(\Omega))}&\leq \liminf\limits_{m\to\infty}\norm*{g_m}_{L^\infty(0,t;L^2(\Omega))}, \label{eq:pqm_tool}\\
\norm*{g}_{L^2(0,t;L^2(\Omega))}&\leq \liminf\limits_{m\to\infty}\norm*{g_m}_{L^2(0,t;L^2(\Omega))}. \nonumber
\end{align}
Using the converges results in \eqref{eq:pqmconvtable} and \eqref{eq:pqmcrossreallim}, as in \cite[Subsection 8.6]{Brunk.d}, we can pass to the limit and obtain
\begin{align}
E(t)=&\left(\int_\Omega \frac{c_0}{2}\snorm{\nabla\phi(t)}^2 + F(\phi(t)) + \frac{1}{2}|q(t)|^2 + \frac{1}{2}|\u(t)|^2 \dx\right) \nonumber \\
& +\iQt \snorm*{n(\phi)\na\mu - \na\big(A(\phi)q\big)}^{2}  +\varepsilon_1\snorm*{\nabla q}^2+  \frac{1}{\tau(\phi)}q^2 + \snorm{\eta^{1/2}(\phi)\Du}^2\dx\ddta   \\
\leq&  \left(\int_\Omega \frac{c_0}{2}\snorm{\nabla\phi(0)}^2 + F(\phi(0)) + \frac{1}{2}|q(0)|^2 + \frac{1}{2}|\u(0)|^2 \dx\right). \nonumber
\end{align}

\section{Weak-Strong Uniqueness}
In this section, we state our result, that the viscoelastic phase separation model \eqref{eq:full_model}  has a relative energy structure, which implies weak-strong uniqueness.

In order to prove our result, we need some further assumptions on the parameter functions.

\begin{Assumptions}
	\label{ass:weakstrong}
	We assume the following additional regularity for the parametric functions and the potential function:
	\begin{itemize}
		\item $n(s), \eta(s), \tau(s)$ are $C^1(\mathbb{R})$ functions such that $\norm{n'}_{L^\infty(\mathbb{R})},\norm{\eta'}_{L^\infty(\mathbb{R})},\norm{\tau'}_{L^\infty(\mathbb{R})}\leq c$.
		\item $A(s)$ is a $C^2(\mathbb{R})$ function with $\norm*{A^\prime}_{L^\infty(\mathbb{R})}\leq c$ and $\norm*{A^{\prime\prime}}_{L^\infty(\mathbb{R})}\leq c$.
		\item We assume that $F\in C^3(\mathbb{R})$ and $\snorm*{F^{\prime\prime\prime}(x)}\leq c_{1,3}\snorm*{x}+c_{2,3},\quad c_{1,3},c_{2,3}\geq 0$.
  \item We introduce the notation $m(s)=n^2(s)$.
	\end{itemize} 
\end{Assumptions}

In principle, all relevant calculations are already done in \cite{brunk2021relative}. We have to verify several regularity assumptions for the weak solution and compute the necessary regularity of the strong solutions.

By virtue of the weak formulation \eqref{eq:weak_sol_full} we can see that the equation for the chemical potential $\mu$, cf. $\eqref{eq:weak_sol_full}_3$ can be rewritten as
\begin{equation}
\int_\Omega \(\mu-F^\prime(\phi)\)\xi \dx = -c_0\int_\Omega \Delta\phi\xi \dx, \text{ for all } \xi\in H^1(\Omega).
\end{equation}
Since we know that $\na\mu$ is bounded in $L^{5/3}(\Omega_T)$, cf. \eqref{eq:weak_sol_full}, we obtain 
\begin{align}
\int_0^T\int_\Omega \snorm*{F^{\prime\prime}(\phi)\na\phi}^{5/3} \dx\dt &\leq c \int_0^T\int_\Omega \snorm*{(1+\snorm*{\phi}^2)\snorm*{\na\phi}}^{5/3} \dx\dt \nonumber\\
&\leq c\norm*{\na\phi}_{L^{5/3}(L^{5/3})}^{5/3} + c\int_0^T \norm*{\phi}_{6}^{10/3}\norm*{\na\phi}_{15/4}^{5/3}\dt \nonumber\\
& \leq c\norm*{\na\phi}_{L^{5/3}(L^{5/3})}^{5/3} + c\norm*{\phi}_{L^\infty(H^1)}^{10/3}\norm*{\na\phi}_{L^{5/3}(L^4)}^{5/3}. \label{eq:relenw3_reg}
\end{align}
The second inequality is obtained by using Hölder inequality with exponents $9/5$ and $9/4$.
Due to the regularity of the weak solution \eqref{eq:weak_sol_full} we observe that (\ref{eq:relenw3_reg}) is bounded, and therefore we find that $\na F^{\prime}(\phi)$ is bounded in $L^{5/3}(\Omega_T)$. Using the elliptic regularity we find
\begin{equation}
\phi \in L^{5/3}(0,T;W^{3,5/3}(\Omega)).
\end{equation}

\begin{Definition}
	\label{defn:more_reg_sol}
	A quadruple $(\psi,\pi,Q,\U)$ is called a \emph{more regular weak solution} if it is a weak solution in the sense of Definition \ref{defn:weak_sol_full} and the following holds
	\begin{align}
		\psi &\in L^2(0,T;H^3(\Omega))\cap W^{1,5/2}(0,T;(W^{1,5/3}(\Omega))^*),\nonumber \\
		Q &\in L^4(0,T;L^\infty(\Omega))\cap L^4(0,T;H^1(\Omega))\cap L^4(0,T;W^{1,6}(\Omega))\cap H^1(0,T;(H^{1}(\Omega))^*),\nonumber\\
		\U &\in L^4(0,T;L^\infty(\Omega)^3)\cap L^4(0,T;V)\cap L^2(0,T;W^{1,3}(\Omega)^3)\cap H^1(0,T;V^*),
	\end{align}
	\begin{align}
		\int_0^t &
		\norm*{n(\psi)\na\pi - \na\big(A(\psi)Q \big) }_4^4 +\norm*{\mathrm{div}\big[n(\psi)\na\pi - \na\big(A(\psi)Q\big)\big]}_2^2 \nonumber \\
		&+ \norm*{\div{\U\psi-n^2(\psi)\na\pi+n(\psi)\na\big(A(\psi)Q\big)}}_2^2 \leq C. \label{eq:reg_strong}
	\end{align}
\end{Definition}

In analogy to \cite[Definition 4.2]{brunk2021relative} we state the following lemma.
\begin{Lemma}\label{lem:timeid}
	Let $(\phi,\mu,q,\u)$ be a global weak solution of \eqref{eq:weak_sol_full} in the sense of Definition \ref{defn:weak_sol_full} and $(\psi,Q,\U)$ functions from the following spaces
	\begin{align*}
	\psi &\in L^2(0,T;H^3(\Omega)),& \quad Q &\in L^{4}(0,T;H^{1}(\Omega)),& \quad \U &\in L^4(0,T;V), \\
	\pt\psi &\in L^{5/2}(0,T;(W^{1,5/3}(\Omega))^*),& \quad \pt Q &\in L^2(0,T;(H^{1}(\Omega))^*),& \quad \pt\U &\in L^2(0,T;V^*).
	\end{align*}
	The following identities hold for a.a. $t,s\in(0,T)$	
	\begin{align*}
	\int_\Omega \na\phi(t)\cdot\na\psi(t) - \na\phi(s)\cdot\na\psi(s) \dx &= -\int_s^t \int_\Omega \pt\phi(\dta)\Delta\psi(\dta) + \Delta\phi(\dta)\pt\psi(\dta) \dx\,\mathrm{d}\dta, \\
	\int_\Omega \phi(t)\psi(t) - \phi(s)\psi(s) \dx &= \int_s^t \int_\Omega \pt \phi(\dta)\psi(\dta) + \phi(\dta)\pt\psi(\dta) \dx\,\mathrm{d}\dta, \\
	\int_\Omega q(t) Q(t) - q(s)Q(s) \dx &= \int_s^t \int_\Omega \pt q(\dta)Q(\dta) + q(\dta)\pt Q(\dta) \dx\,\mathrm{d}\dta, \\
	\int_\Omega \u(t)\cdot\U(t) - \u(s)\cdot\U(s) \dx &= \int_s^t \int_\Omega \pt\u(\dta)\cdot\U(\dta) + \u(\dta)\cdot\pt\U(\dta) \dx\,\mathrm{d}\dta.
	\end{align*}
\end{Lemma}	
The regularity assumptions are chosen such that the integrals involving the weak and more regular solutions make sense.
Then we define a relative energy by
	\begin{align}
	&\mathcal{E}(\phi,q,\u\vert \psi,Q,\U) = \mathcal{E}_{mix}(\phi\vert\psi) + \mathcal{E}_{bulk}(q\vert Q) + \mathcal{E}_{kin}(\u\vert\U) , \label{eq:relen}\\
	&\mathcal{E}_{mix}(\phi\vert\psi) = \int_\Omega \frac{c_0}{2}\snorm*{\na\phi-\na\psi}^2 + F(\phi) - F(\psi) - F^\prime(\psi)(\phi-\psi) + a(\phi-\psi)^2 \dx, \nonumber\\
	&\mathcal{E}_{bulk}(q\vert Q) = \int_\Omega \frac{1}{2}\snorm*{q-Q}^2 \dx, \quad \mathcal{E}_{kin}(\u\vert\U) = \int_\Omega \frac{1}{2}\snorm*{\u-\U}^2 \dx. \nonumber
	\end{align}
One can see, that for $a>c_4/2$ the relative energy can be bounded from below by the natural energy norms of \eqref{eq:weak_sol_full}.

\begin{tcolorbox}
\begin{Theorem}[Relative energy and weak-strong uniqueness, {\cite[Theorem 4.5]{brunk2021relative}}]
	\label{theo:wsu_reg}
	Let Assumptions \ref{ass:paramreg} and \ref{ass:weakstrong} hold. Let $(\phi,\mu,q,\u)$ be a weak solution corresponding to the initial data $(\phi_0,q_0,\u_0)$ in the sense of Definition \ref{eq:weak_sol_full}. Let $(\psi,\pi,Q,\U)$ be a more regular solution in the sense of Definition \ref{defn:more_reg_sol} starting from the initial data $(\psi_0,Q_0,\U_0)$ and exists up to time $T^\dagger\leq T$. 
	Then the following inequality holds
	\begin{align}
	\mathcal{E}(t) + \tfrac{1}{2}\mathcal{D} \leq \mathcal{E}(0)  + c\int_0^t g(\tau)\mathcal{E}(\tau)\ddta, \label{eq:re1}
	\end{align}	
	with  $g\in L^1(0,T^\dagger)$ and $\mathcal{D}$ is given by
	\begin{align*}
	\mathcal{D} =& \iQt \eta(\phi)\snorm*{\Du - \mathrm{D}\U}^2 + \frac{1}{\tau(\phi)}(q-Q)^2+\varepsilon_1\snorm*{\na q - \na Q}^2\dx\ddta \\
	&+ \iQt  \Big|n(\phi)(\na\mu-\na\pi) - \na\big(A(\phi)(q-Q) \big)\Big|^2 \dx\ddta. \nonumber
	\end{align*}
	
	Furthermore, if $(\phi_0,q_0,\u_0)=(\psi_0,Q_0,\U_0)$ holds a.e. in $\Omega$. Then $(\phi,\mu,q,\u)=(\psi,\pi,Q,\U)$ a.e. in $\Omega\times[0,T^\dagger]$.
\end{Theorem}
\end{tcolorbox}
\begin{proof}
The result and proof are given in \cite[Theorem 4.5]{brunk2021relative}.
\end{proof}

\section{Degenerate case}
In this section, we want to focus on the degenerate case in three space dimensions. First, we will provide an existence result and afterwards we will discuss the weak-strong uniqueness principle.  The structure of the existence proof is analogous to the two-dimensional case, see \cite{Brunk.c}. However, we will repeat the relevant ideas and mainly focus on the regularity loss and the limiting process.

First, we state a set of assumptions for the degenerate case.
\begin{Assumptions} \label{ass:deg}
	\hspace{1em} \\
	\vspace{-1em}
	\begin{itemize}
        \item We assume that $n\in C^1([0,1])$ with $n(x)=0$ if and only if $x\in\{0,1\}$. The mobility functions $n(x), n^2(x)$ are continuously extended by zero on $\mathbb{R}\setminus[0,1]$.
		\item The potential can be divided into $F=F_1+F_2$ with a convex part $F_1\in C^2(0,1)$ and a concave part $F_2\in C^2([0,1])$. $F_2$ is continuously extended on $\mathbb{R}$ such that $\norm*{F_2^{\prime\prime}}_{L^\infty(\R)}\leq F_0$.
		\item The convex part $F_1$ additionally satisfies $(n^2F_1^{\prime\prime})\in C([0,1])$.
		\item We assume $\norm*{\frac{A(x)}{n(x)}}_{L^\infty(\mathbb{R})}\leq c,\quad \norm*{\frac{A^\prime(x)}{n(x)}}_{L^\infty(\mathbb{R})}\leq c.$
		\item The boundary condition $\restr{n^2(\phi)\na\mu\cdot\mathbf{n}}{\p\Omega}=0$ holds.
	\end{itemize}
\end{Assumptions}

The assumptions on $A, A'$ are understood in the following sense. Since both functions are bounded by Assumptions \ref{ass:paramreg} this forces the function $A, A'$ to decay in a similar or even faster way than $n$. From a physical point of view, this is acceptable, since on the pure phases $\phi \pm 1$ the Cahn-Hilliard type equation should reduce to simple transport.
For the convenience of the upcoming analysis, we set $m=n^2$.

\begin{Definition}
	\label{defn:weak_sol_deg_full}	
	Let the initial data $(\phi_0,q_0,\u_0)\in H^1(\Omega)\times L^2(\Omega)\times H.$ Then for every $T>0$ the quadruple $(\phi,q,\J,\u)$ is called a weak solution of (\ref{eq:full_model}) if it satisfies
	\begin{align*}	
	&\phi \in L^{\infty}(0,T;H^1(\Omega))\cap L^2(0,T;H^2(\Omega))\cap H^1(0,T;(H^{1}(\Omega))^*),\\
	&q\in L^\infty(0,T;L^2(\Omega))\cap L^2(0,T;H^1(\Omega))\cap W^{1,4/3}(0,T;(H^{1}(\Omega))^*), \\
	&\u\in L^{\infty}(0,T;H)\cap L^2(0,T;V)\cap W^{1,4/3}(0,T;V^*),\\
	&\J=n(\phi)\Jh,\Jh\in L^{5/3}(0,T;L^{5/3}(\Omega)),\quad A(\phi)q\in L^{5/3}(0,T;W^{1,5/3}(\Omega)).
	\end{align*}
	Further, for any test function $(\psi,\zeta,\bxi,\vv)\in H^1(\Omega)\times H^{1}(\Omega)\times H^1(\Omega)\cap L^\infty(\Omega)\times V$ and almost every $t\in(0,T)$ it holds
	\begin{align}
	&\int_\Omega\pt\phi\psi \dx + \int_\Omega\(\u\cdot\nabla\phi\)\psi \dx + \int_\Omega n(\phi)\(\Jh - \nabla\big(A(\phi)q\big) \)\cdot\nabla\psi \dx =0 \nonumber\\
	&\int_\Omega\pt q\zeta \dx + \int_\Omega\(\u\cdot\nabla q\)\zeta \dx + \int_\Omega\frac{q\zeta}{\tau(\phi)}\dx  + \int_\Omega \varepsilon_1\nabla q\cdot\nabla\zeta\dx \nonumber \\
	&\hspace{12.75em}+  \int_\Omega\Big[\nabla\big(A(\phi)q\big) - \Jh\big]\cdot\na\big(A(\phi)\zeta\big)\dx= 0   \nonumber\\
	&\int_\Omega \J\cdot\bxi \dx = c_0\int_\Omega \Delta\phi\div{n^2(\phi)\bxi}\dx + \int_\Omega n^2(\phi)F^{\prime\prime}(\phi)\nabla\phi\cdot\bxi \dx\label{eq:weak_sol_sys} \\
	&\int_\Omega\pt\u\cdot\vv \dx + \int (\u\cdot\na)\u\cdot\vv \dx - \int_\Omega\eta(\phi)\Du:\mathrm{D}\vv \dx  + \int_\Omega c_0\Delta\phi\nabla\phi\cdot\vv \dx = 0 \nonumber
	\end{align}
	and the initial data, i.e. $(\phi(0),q(0),\u(0))=(\phi_0,q_0,\u_0)$ are attained.
\end{Definition}
The above formulation for $\J$ is a weak version of $\J = n^2(\phi)\big(-c_0\nabla\Delta\phi +F^{\prime\prime}(\phi)\nabla\phi\big)$ thus for a smooth solution we can identify $\J=n^2(\phi)\nabla\mu$. \\

\begin{tcolorbox}
	\begin{Theorem}\label{theo:deg1}
		Let Assumptions \ref{ass:deg}  hold. Further, let $\phi_0:\Omega\to [0,1]$ and $\phi_0\in H^1(\Omega)$. The potential function $F$ and the entropy function $G$, cf. (\ref{eq:defentropy}), fulfil
		\begin{equation}
		\int_\Omega \Big( F(\phi_0) + G(\phi_0) \dx \Big) < +\infty. \label{eq:startres}
		\end{equation}
		Then for any given $T < \infty$ there exists a global weak solution $(\phi,q,\J,\u)$ of the viscoelastic phase separation model (\ref{eq:full_model}) in the sense of Definition \ref{defn:weak_sol_deg_full}. Moreover,
		\begin{itemize}
			\item the integrated energy inequality 
			\begin{align}
			& \(\int_\Omega \frac{c_0}{2}\snorm*{\na\phi(t)}^2 + F(\phi(t)) + \frac{1}{2}\snorm*{q(t)}^2 + \frac{1}{2}\snorm*{\u(t)}^2  \dx \) \label{eq:fullie1_deg} \\
			&+ \iQt \snorm*{\Jh - \na\big(A(\phi)q\big) }^2 + \frac{1}{\tau(\phi)}|q|^2 + \varepsilon_1\snorm*{\na q}^2 + \eta(\phi)\snorm*{\Du}^2\dx\ddta\nonumber \\
			&\leq  \(\int_\Omega \frac{c_0}{2}\snorm*{\na\phi(0)}^2 + F(\phi(0))+ \frac{1}{2}\snorm*{q(0)}^2 + \frac{1}{2}\snorm*{\u(0)}^2 \dx \). \nonumber
			\end{align}	
			 holds.
			\item $\phi(x,t) \in [0,1] \text{  for a.e. } (x,t)\in \Omega\times(0,T). $
		\end{itemize}
		If the mobility function satisfies $(n^2)^\prime(0) = (n^2)^\prime(1) = 0$, then for a.e. $t\in (0,T)$ the set
		\begin{equation*}
		\{(x,t)\in\Omega\times (0,T) \mid \phi(x,t) = 0 \text{ or } \phi(x,t) = 1\}
		\end{equation*}
		has zero measure.
	\end{Theorem}
\end{tcolorbox}

The proof is structured as follows, cf. \cite[Proof of Theorem 4.1]{Brunk.c}:
\begin{enumerate}
    \item Existence of a sequence of weak solutions by regularizing the mobilities $m,n$, the potential $F$ and the bulk modulus A, such that Theorem \ref{theo:ex_full} holds.
    \item A priori estimate using the energy and entropy estimates.
    \item $L^\infty$-bounds for $\phi$, by exploiting a priori bounds and the singular behaviour of entropy $G$.
    \item Passage to the limit in the weak formulation and the energy.
    \item Extended bounds for $\phi$ for stronger degenerating mobilities, using again the singular behaviour of $G$.
\end{enumerate}

\subsection{Proof}

\noindent\textbf{Regularization:}\\
We start by introducing a suitable regularized problem and approximate the degenerate mobility $m$, the logarithmic potential $F$ and the bulk modulus by a non-degenerate mobility $m_\delta$, a smooth potential $F_\delta$ and $A_\delta$ with a parameter $\delta\in (0,\frac{1}{2})$.
\begin{equation}
\label{eq:appmob}
m_\delta(s) = \left\{\begin{array}{ll}
m(\delta), &\text{if } s\leq \delta\\
m(s), &\text{if } \delta \leq s \leq 1-\delta \\
m(1-\delta), &\text{if } s \geq 1- \delta.
\end{array}\right.
\end{equation}
The regularization for $n_\delta$ is similar and $n_\delta^2=m_\delta$.
Since $F_2$ is already defined on $\mathbb{R}$ and bounded, cf. Assumptions \ref{ass:deg}, we set $F_{2,\delta}=F_2$. Further,
\begin{align}
\label{eq:apppot}
F_{1,\delta}(1/2) = F_1(1/2), \qquad F_{1,\delta}^\prime (1/2) = F^\prime_1(1/2), \text{ and } \\
\label{eq:apppot2}
F_{1,\delta}^{\prime\prime}(s) = \left\{\begin{array}{ll}
F_1^{\prime\prime}(\delta), &\text{if } s\leq \delta\\
F_1^{\prime\prime}(s), &\text{if } \delta \leq s \leq 1-\delta \\
F_1^{\prime\prime}(1-\delta), &\text{if } s \geq 1- \delta.
\end{array}\right.
\end{align}
We note that $F_1(s) = F_{1,\delta}(s)$ for $s\in[\delta,1-\delta]$. The regularized bulk modulus $A_\delta$ is defined, such that $A_\delta/n_\delta$ and $A'_\delta/n_\delta$ are bounded in $L^\infty(\mathbb{R})$. Thus, due to Assumptions \ref{ass:paramreg} the limit is also bounded.\\

The system \eqref{eq:full_model} with $F_\delta$ and $m_\delta$ fulfils the hypothesis of Theorem \ref{theo:ex_full} for every $\delta\in(0,\frac{1}{2})$. Consequently, there is a sequence of regularized weak solutions, denoted by $(\phi_{\delta},q_{\delta},\mu_{\delta},\u_{\delta})$, such that
\begin{align*}
&\phi_{\delta} \in L^{\infty}(0,T;H^1(\Omega))\cap L^2(0,T;H^2(\Omega))\cap H^1(0,T;(H^{1}(\Omega))^*),\\
& q_{\delta}\in L^\infty(0,T;L^2(\Omega))\cap L^2(0,T;H^1(\Omega))\cap W^{1,4/3}(0,T;(H^{1}(\Omega))^*),\\
&\u_{\delta}\in L^{\infty}(0,T;H)\cap L^2(0,T;V)\cap W^{1,4/3}(0,T;V^*),\\
&A(\phi_\delta)q_\delta\in L^{5/3}(0,T;W^{1,5/3}(\Omega)), \quad
J_\delta:=m_\delta(\phi_\delta)\na\mu_\delta \in L^{5/3}(0,T;L^{5/3}(\Omega)).
\end{align*}
The family $(\phi_\delta,q_\delta,\u_\delta)$ fulfils the weak formulation (\ref{eq:weak_sol_full}) for every $\delta > 0$, which in the following will be denoted by $(\ref{eq:weak_sol_full})_\delta$. Therefore, the above estimates are independent of $\delta$. However, due to the degeneracy of $m_\delta$, we will lose control of $\na\mu_\delta.$
The key step is now to obtain $m(\phi)$ independent estimates. To this end, we construct the so-called entropy function \cite{Boyer.1999,Elliott.1996} $G$ via
\begin{equation}
G_\delta(1/2) = 0, \quad G_\delta^\prime(1/2) = 0 \quad G_\delta^{\prime\prime}(s) = \frac{1}{m_\delta(s)}, \text{ for all } s\in\mathbb{R}. \label{eq:defentropy}
\end{equation}
Following the calculations presented for the regular case, see also \cite{Brunk.c},  we derive the estimates
\begin{align*}
\norm*{\Delta\phi_\delta}_{L^2(L^2)} + \norm*{\int_\Omega G_\delta(\phi_\delta) \dx}_{L^\infty(0,T)} \leq c.
\end{align*} 

\noindent\textbf{$L^\infty$-bounds for $\phi$:} \\
The following Lemma from \cite[Lemma 8.2]{Brunk.c} can be proven with the above bounds.
\begin{Lemma}
	Let $\phi_\delta$ be the solution of the regularized problem $(\ref{eq:weak_sol_full})^\delta_1$. Then it holds that $\phi_\delta(x,t)$ converges uniformly for a.e. $(x,t)$ in $\Omega\times(0,T)$ to $\phi(x,t)\in[0,1]$.
\end{Lemma}
Technically, the lemma uses a priori bounds on $G_\delta$ and its limit, combined with the possible singular behaviour of $G$ at $0$ and $1$.
The next step is the passage to the limit. Following \cite{Brunk.c} we observe that the passage is also valid, and we will only comment on the crucial points.

The integral 
\begin{equation*}
\int_0^T \int_\Omega (\J_\delta-\J)\cdot\bxi\varphi(t) \dx \ddta    
\end{equation*}
goes to zero due to the weak convergence of $\J_\delta$ in $L^{5/3}(\Omega_T)$, since the test function $\bxi\varphi(t)\in L^\infty(0,T;H^1(\Omega)\cap L^\infty(\Omega))$.

The crucial integrals for the bulk stress equation converge, since $A_\delta(\phi_\delta), A'_\delta(\phi_\delta)$ converge almost everywhere. Hence, the integrals converge similarly as in the regular case.

Finally, we observe that the capillary stress tensor, connecting the Cahn-Hilliard equation with the Navier-Stokes equations, converges more easily than in the regular case since $\mu_\delta$ is not present anymore.

The limit in the energy inequality follows from the same arguments as in the regular case.

\subsection{Conditional weak-strong uniqueness}
In contrast to Theorem \ref{theo:wsu_reg}, the degenerate weak solution is not regular enough to apply the relative energy structure. Even if the weak solution would be regular enough, the bounds for the relative energy in Theorem \ref{theo:wsu_reg} depend inversely on the lower bounds of the mobility. This blows up in the degenerate case. Therefore, we have to assume additional properties of the weak solution. Indeed, a natural assumption is the strong separation principle, i.e.  $\phi$ satisfies
\begin{equation}
\phi\in[\kappa,1-\kappa] \text{ a.e. in } \Omega \text{ for almost all } t\in(0,T). \label{eq:sep}
\end{equation}
Here the constant $\kappa>0$ will depend on time and the data.

\begin{tcolorbox}
	\begin{Theorem}
		\label{theo:wsu_deg}
		Let Assumptions \ref{ass:deg}, \ref{ass:weakstrong} hold. Let $(\phi,\mu,q,\u)$ be a weak solution starting from the initial data $(\phi_0,q_0,\u_0)$ in the sense of Definition \ref{defn:weak_sol_deg_full}. Let $(\psi,\pi,Q,\U)$ be a more regular solution in the sense of Definition \ref{defn:more_reg_sol} starting from the initial data $(\psi_0,Q_0,\U_0)$ existing up to time $T^\dagger\leq T$. Let $\phi$ and $\psi$ satisfy the strict separation principle \eqref{eq:sep}. Then the relative energy inequality \eqref{eq:re1} holds and for coinciding initial data we have $(\phi,\mu,q,\u)=(\psi,\pi,Q,\U)$ a.e. in $\Omega\times[0,T^\dagger]$. 
	\end{Theorem}
\end{tcolorbox}
\begin{Remark}\phantom{e}
\begin{itemize}
		\item A rigorous proof of the separation principle for Cahn-Hilliard with degenerate mobility and logarithmic potential is still open. For constant mobilities and sufficiently singular potentials, the result was proven in \cite{Cherfils.2011,Miranville.2004}.  
		\item The separation principle might be obtained by refining the higher entropy estimate in \cite{Gruenentropy} to the mobility functions used here. However, this is part of further research.
		\end{itemize}
\end{Remark}

\begin{proof}
Reviewing the last section, the strict separation principle \eqref{eq:sep} implies that $n(\phi)^{-1}$, $m(\phi)^{-1}$, $F'(\phi)\in L^\infty(\Omega_T)$. This immediately implies the solution is non-degenerate and fulfils the weak formulation \eqref{eq:weak_sol_full}, see for instance \cite{Dai.2016}. In detail, the above bound implies that
\begin{equation*}
L^p(\Omega_T) \ni \Jh = n(\phi)\na\mu \Longrightarrow \norm*{\na\mu}_{L^p(L^p)}\leq c.
\end{equation*}
At this point, we can apply Theorem \ref{theo:wsu_reg}. The lower bound of the mobility is replaced by the lower bound restricted to the separation interval in \eqref{eq:sep}.
\end{proof}

\section{Conclusions}
In this work, we have proven the existence of global-in-time weak solutions for the viscoelastic phase separation model \eqref{eq:full_model} in three space dimensions. The main difficulty is dealing with the lost regularity due to the highly nonlinear structure, in contrast to the two-dimensional case. Furthermore, we apply our concept of relative energy developed on \cite{brunk2021relative} in three space dimensions, which yields the weak-strong uniqueness principle.  Finally, we have proven the existence of global weak solutions in the degenerate case, which generalizes \cite{Brunk.c} to three-space dimensions. Moreover, we have discussed that by requiring the strict separation principle for the weak solution, the relative energy method can be applied. This implies a conditional weak-strong-uniqueness result. A rigorous proof of this separation principle in the case of degenerate mobilities is still open.

\section*{Acknowledgment}
Funded by the Deutsche Forschungsgemeinschaft (DFG, German Research Foundation) - Project number 233630050 - TRR 146.

\bibliography{paper_bib}

\begin{thebibliography}{10}

\bibitem{Barrett}
J.~W. Barrett and E.~S\"{u}li.
\newblock Existence and equilibration of global weak solutions to kinetic
  models for dilute polymers {I}: {Finitely} extensible nonlinear bead-spring
  chains.
\newblock {\em Math Models Methods Appl Sci}, 21(06):1211--1289, 2011.

\bibitem{MichalBathory.2020}
M.~Bathory, M.~Bulíček, and J.~Málek.
\newblock Large data existence theory for three-dimensional unsteady flows of
  rate-type viscoelastic fluids with stress diffusion.
\newblock {\em Adv Nonlinear Anal}, 10(1):501--521, 2021.

\bibitem{Boyer.1999}
F.~Boyer.
\newblock {Mathematical study of multiphase flow under shear through order
  parameter formulation}.
\newblock {\em {Asymptotic Anal}}, 20:175--212, 1999.

\bibitem{Brunk2021}
A.~Brunk, B.~D\"{u}nweg, H.~Egger, O.~Habrich,
  M.~Luk{\'a}{\v{c}}ov{\'a}-Medvi{\softd}ov{\'a}, and D.~Spiller.
\newblock Analysis of a viscoelastic phase separation model.
\newblock {\em J Condens Matter Phys}, 33(23):234002, 2021.

\bibitem{Brunk.d}
A.~Brunk and M.~Luk{\'a}{\v{c}}ov{\'a}-Medvi{\softd}ov{\'a}.
\newblock {Global existence of weak solutions to the two-phase viscoelastic
  phase separation: Part I Regular Case}.
\newblock {\em Nonlinearity}, 35(7):3417, 2022.

\bibitem{Brunk.c}
A.~Brunk and M.~Luk{\'a}{\v{c}}ov{\'a}-Medvi{\softd}ov{\'a}.
\newblock {Global existence of weak solutions to the two-phase viscoelastic
  phase separation: Part II Degenerate Case}.
\newblock {\em Nonlinearity}, 35(7):3459, 2022.

\bibitem{brunk2021relative}
A.~Brunk and M.~Lukáčová-Medvid'ová.
\newblock Relative energy and weak–strong uniqueness of a two-phase
  viscoelastic phase separation model.
\newblock {\em Z Angew Math Mech}, page e202100240, 2022.

\bibitem{Cherfils.2011}
L.~Cherfils, A.~Miranville, and S.~Zelik.
\newblock {The Cahn-Hilliard equation with logarithmic potentials}.
\newblock {\em {Milan J Math}}, 79(2):561--596, 2011.

\bibitem{Dai.2016}
S.~Dai and Q.~Du.
\newblock {Weak Solutions for the Cahn--Hilliard Equation with Degenerate
  Mobility}.
\newblock {\em {Arch Ration Mech An}}, 219(3):1161--1184, 2016.

\bibitem{Gruenentropy}
R.~Dal~Passo, H.~Garcke, and G.~Gr\"{u}n.
\newblock On a fourth-order degenerate parabolic equation: Global entropy
  estimates, existence, and qualitative behavior of solutions.
\newblock {\em SIAM J. Math. Anal.}, 29(2):321–342, Mar. 1998.

\bibitem{Groot.2016}
S.~R. de~Groot and P.~Mazur.
\newblock {\em {Non-equilibrium thermodynamics}}.
\newblock {Dover Books on Physics}. {Dover Publications, Inc}, New York, 2016.

\bibitem{Elliott.1996}
C.~M. Elliott and H.~Garcke.
\newblock {On the Cahn--Hilliard equation with degenerate mobility}.
\newblock {\em {SIAM J Math Anal}}, 27(2):404--423, 1996.

\bibitem{Feireisl.2016}
E.~Feireisl.
\newblock {On weak solutions to a diffuse interface model of a binary mixture
  of compressible fluids}.
\newblock {\em {DCDS-S}}, 9(1):173--183, 2016.

\bibitem{Feireisl.2018}
E.~Feireisl, Y.~Lu, and A.~Novotn{\'y}.
\newblock {Weak-strong uniqueness for the compressible Navier-Stokes equations
  with a hard-sphere pressure law}.
\newblock {\em {Sci China Math}}, 61(11):2003--2016, 2018.

\bibitem{Feireisl.2012}
E.~Feireisl and A.~Novotn{\'y}.
\newblock {Weak--Strong uniqueness property for the full
  Navier--Stokes--Fourier system}.
\newblock {\em {Arch Ration Mech An}}, 204(2):683--706, 2012.

\bibitem{Feireisl.2019}
E.~Feireisl, M.~Petcu, and D.~Pra{\v{z}}{\'a}k.
\newblock {Relative energy approach to a diffuse interface model of a
  compressible two--phase flow}.
\newblock {\em {Math Method Appl Sci}}, 42(5):1465--1479, 2019.

\bibitem{feireisl2021numerical}
E.~Feireisl, M.~Petcu, and B.~She.
\newblock Numerical analysis of a model of two phase compressible fluid flow,
  2021.

\bibitem{Giesselmann2017}
J.~Giesselmann and A.~E. Tzavaras.
\newblock Stability properties of the {E}uler{\textendash}{K}orteweg system
  with nonmonotone pressures.
\newblock {\em Appl Anal}, 96(9):1528--1546, Jan. 2017.

\bibitem{Grmela.1997b}
M.~Grmela and H.~C. {\"O}ttinger.
\newblock {Dynamics and thermodynamics of complex fluids. I. Development of a
  general formalism}.
\newblock {\em {Phys Rev E}}, 56(6):6620--6632, 1997.

\bibitem{Grmela.1997}
M.~Grmela and H.~C. {\"O}ttinger.
\newblock {Dynamics and thermodynamics of complex fluids. II. Illustrations of
  a general formalism}.
\newblock {\em {Phys Rev E}}, 56(6):6633--6655, 1997.

\bibitem{Grun.1995}
G.~Gr{\"u}n.
\newblock {Degenerate parabolic differential equations of fourth order and a
  plasticity model with non-local hardening}.
\newblock {\em {Z Anal Anwend}}, 14(3):541--574, 1995.

\bibitem{Hohenberg.1977}
P.~C. Hohenberg and B.~I. Halperin.
\newblock {Theory of dynamic critical phenomena}.
\newblock {\em {Rev Mod Phys}}, 49(3):435--479, 1977.

\bibitem{Jungel.2016}
A.~J{\"u}ngel.
\newblock {\em {Entropy methods for diffusive partial differential equations}}.
\newblock {BCAM Springer Briefs}. Springer, 2016.

\bibitem{Miranville.2004}
A.~Miranville and S.~Zelik.
\newblock {Robust exponential attractors for Cahn-Hilliard type equations with
  singular potentials}.
\newblock {\em {Math Meth Appl Sci}}, 27(5):545--582, 2004.

\bibitem{Tanaka.200}
H.~Tanaka.
\newblock {Viscoelastic phase separation}.
\newblock {\em {J. Phys.: Condens. Matter}}, 12(15):R207, 200.

\bibitem{Zhou.2006}
D.~Zhou, P.~Zhang, and W.~E.
\newblock {Modified models of polymer phase separation}.
\newblock {\em {Phys Rev E}}, 73(6):061801, 2006.

\end{thebibliography}
\bibliographystyle{abbrv}

\end{document}